\documentclass{proc-l}

\usepackage{epsf}
\usepackage{epsfig}
\usepackage{graphicx}
\usepackage{xcolor}
\usepackage{bm}
\usepackage{tabularx}
\usepackage{amstext}
\usepackage{float}
\usepackage{array, nicefrac, caption}

\newcommand{\1}{{\bf 1}}
\newcommand{\divprop}{\lfloor}

\newtheorem{theorem}{Theorem}

\newtheorem{conjecture}{Conjecture}

\allowdisplaybreaks
\begin{document}

\frenchspacing

\title{Number of ordered factorizations and recursive divisors}

\author{T. M. A. Fink}
\address{London Institute for Mathematical Sciences, Royal Institution, 21 Albermarle St, London W1S 4BS, UK}

\date{\today}
\vspace{-0.2in}
\begin{abstract}
The number of ordered factorizations and the number of recursive divisors are two related arithmetic functions that are recursively defined. 
But it is hard to construct explicit representations of these functions. 
Taking advantage of their recursive definition and a geometric interpretation, we derive three closed-form expressions for them both.
These expressions shed light on the structure of these functions and their number-theoretic properties.
Surprisingly, both functions can be expressed as simple generalized hypergeometric functions.
\vspace*{-0.25in}
\end{abstract}
\maketitle
\section{Introduction}
\noindent
This paper is devoted to two related arithmetic functions that are recursively defined.
The first is the number of ordered factorizations into integers greater than one.
The second is the number of recursive divisors, 
which measures the extent to which a number $n$ is highly divisible, whose quotients are highly divisible, and so on.
\\ \indent
The first function was introduced 90 years ago by Kalmar \cite{Kalmar31}, and for this reason is called $K(n)$.
For example, $K(8) = 4$, since 8 can be factorized in 4 ways: $8 = 2 \cdot 4 =  4 \cdot 2 = 2 \cdot 2 \cdot 2$.
Other values of $K(n)$ are given in Table 1.
Hille \cite{Hille36} extended Kalmar's results and gave them prominence.
Canfield et al. \cite{Canfield83} and Del\'{e}glise et al. \cite{Deleglise08} studied the indices of sequence records of $K$, that is, values of $n$ for which $K(n) > K(m)$ for all $m < n$. 
Newburg and Naor \cite{Newberg93} showed that $K$ arises in computational biology, in the so-called probed partial digest problem,
which prompted Chor et al. \cite{Chor00} to study the upper bound of $K$.
This bound was improved by Klazar and Luca \cite{Klazar07}, who also considered arithmetic properties of the function.
\\ \indent
$K(n)$ can be defined recursively as follows. 
Let $n = m \, d$, where $m > 1$. 
All of the factorizations of $n$ that begin with $m$ can be had by counting the ordered factorizations of $d$,
of which there are $K(d)$.
Thus we obtain the recursion relation $K(n) = \sum_{d\divprop n} K(d)$, 
where $m \divprop n$ means $m \vert n$ and $m < n$.
Along with the initial condition $K(1) = 1$, this completely determines $K(n)$.
This, at least, is how $K$ was originally defined \cite{Hille36}.
But it is much more fruitful to embed the initial condition into the recursion relation itself:
	\begin{equation}
		K(n) = \varepsilon + \sum_{d\divprop n} K(d),
		\label{Kdef}
	\end{equation}
where $\varepsilon = \lfloor 1/n \rfloor = 1,0,0,0,\ldots$. 
As we shall see, this lets us manipulate the defining equation without having to keep track of the corresponding initial condition.
\\ \indent
In contrast with $K(n)$, the second function $\kappa_0(n)$ is much more recent \cite{Fink21, Fink22}. 
It counts the number of recursive divisors: 
	\begin{equation}
		\kappa_0(n) = 1 + \sum_{d\divprop n} \kappa_0(d).
		\label{kappadef}
	\end{equation}
For example, $\kappa_0(8) = 1 + \kappa_0(1) + \kappa_0(2) + \kappa_0(4) = 8$, and other values of $\kappa_0$ are given in Table 1.
$\kappa_0(n)$ is the simplest case of the more general 
$$\kappa_x(n) = n^x + \sum_{d\divprop n} \kappa_x(d),$$ 
which was introduced as a recursive analogue of the usual divisor function $\sigma_x(n)$ \cite{Fink21}.
Like $K(n)$, $\kappa_0(n)$ depends only on the prime signature of $n$---though this is not the case for other values of $x$.
The analogy with $\sigma_x(n)$ motivated the study of recursively perfect numbers ($\kappa_0(n) = n$) and recursively abundant numbers ($\kappa_0(n) > n$) \cite{Fink22}.
\\ \indent
The two sequences $K$ and $\kappa_0$ are intimately related.
In particular, as we showed in \cite{Fink21}, for $n \geq 2$, $\kappa_0(n) = 2 K(n)$.
Furthermore, 
	\begin{equation*}
		\kappa_0(n) = \sum_{d|n} K(d).
	\end{equation*}
\\ \indent
Both $\kappa_0$ and $K$ have a geometric interpretation:
$\kappa_0$ is the number of squares in the divisor tree of $n$, 
whereas $K$ is the number of squares of size 1 in the divisor tree of $n$ (see Fig. 1E).
The divisor tree is constructed as follows. 
Starting with a square of side length $n$ (Fig. 1A), the main arm of the divisor tree is made up of smaller squares with side lengths equal to the proper divisors of $n$ (Fig. 1B).
For each square in the main arm, a secondary arm is made up of squares with side lengths equal to that square's proper divisors (Fig. 1C).
The process is repeated, creating sub-arms off of sub-arms, until the last sub-arms are of size 1 (Fig. 1E).
\\ \indent
A few words on notation.
As we mentioned, $m \divprop n$ means $m \vert n$ and $m < n$, that is, $m$ is a proper divisor of $n$.
We denote the Dirichlet series of an arithmetic function $f$ by $\widetilde {f}$
and the Dirichlet convolution of two arithmetic functions $f$ and $g$ by 
\begin{align*}
		f \star g = \sum_{d|n} f(d) \, g\left(\frac{n}{d}\right).
\end{align*}
\begin{figure}[b!]
\includegraphics[width=1\textwidth]{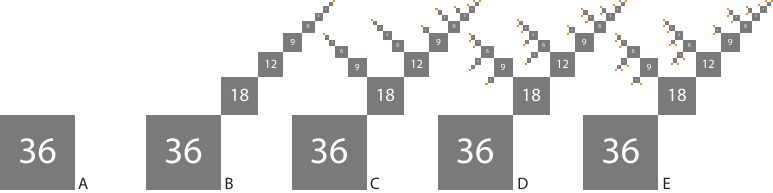} 
\caption{\small
{\bf Divisor trees.}
Both $\kappa_0(n)$ and $K(n)$ have a geometric interpretation.
The number of recursive divisors $\kappa_0(n)$ is the number of squares in the divisor tree of $n$ (E),
whereas the number of recursive divisors $K(n)$ is the number of squares of size 1 (orange squares in E).
As the divisor tree is built up over successive generations, from generation $i = 0$ (A) to $i = \Omega(36) = 4$ (E), 
the number of squares increases by $\upsilon_{i+1}$.
For $n = 36$, the values of $\upsilon_i$ are 
$\upsilon_1 = 1$,
$\upsilon_2 = 8$,
$\upsilon_3 = 19$, 
$\upsilon_4 = 18$ and 
$\upsilon_5 = 6$.
Then $\kappa_0(36) = \upsilon_1 + \ldots + \upsilon_5 = 52$.
As expected, $K(36) = \kappa_0(36)/2 = 26$.
}
\vspace{-14pt}
\label{threetypesplot}
\end{figure}
\section{Statement of results}
\noindent
In this paper, we give three closed form representations of $\kappa_0 = 2 K$ for arbitrary $n$.
As well as making it easier to compute both functions, these representations also shed light on their structure and their number-theoretic properties.
\\ \indent
Let
	\begin{equation*}
		n = p_1^{\alpha_1} p_2^{\alpha_2} \ldots p_\omega^{\alpha_\omega}
	\end{equation*}
be the prime factorization of $n$ and let 
	\begin{equation*}
		\Omega = \alpha_1 + \alpha_2 + \ldots + \alpha_\omega.
	\end{equation*}
\begin{theorem}
	\begin{align*}
		\kappa_0(n) = 2 K(n) &= \frac{1}{2} \sum_{i=0}^\infty \frac{1}{2^i} 
		\prod_{k=1}^\omega \binom{\alpha_k + i}{\alpha_k} \\
				&= \textstyle{\frac{1}{2}} \, {}_\omega F_{\omega-1}\!\left[\begin{matrix}\alpha_1+1, \ldots, \alpha_\omega+1\\1, \ldots, 1\end{matrix} ; \frac{1}{2}\right],
	\end{align*}
\end{theorem}
\noindent
where ${}_\omega F_{\omega-1}$ is the generalized hypergeometric function.
\begin{theorem}
	\begin{align*}
		\kappa_0(n) = 2 K(n) & = \sum_{i=0}^{\Omega} \sum_{j=0}^i (-1)^{i-j} \binom{i}{j} 	
		\prod_{k=1}^\omega \binom{\alpha_k + j}{\alpha_k}.
	\end{align*}
\end{theorem}
\begin{conjecture}
	\begin{align*}
		\kappa_0(n) = 2 K(n) &=	2^{\alpha_\omega} \sum_{i_1=0}^{\alpha_1} 	\ldots 	\sum_{i_{\omega-1}=0}^{\alpha_{\omega-1}} 
		\prod_{k=1}^{\omega-1}	
		\binom{\alpha_k}{i_k} 
		\binom{\alpha_{k+1} + i_1 + \ldots + i_k}{\alpha_{k+1}}.
	\end{align*}
\end{conjecture}
\noindent
We tested the conjecture for $n \in [1,10^5]$ and it is correct in all cases.
\\ \indent
{\bf Example.}
The easiest way to gain some intuition for these three different expressions is to consider an example:
$n = p_1^{\alpha_1} p_2^{\alpha_2} p_3^{\alpha_3} p_4^{\alpha_4}$.
Then $\omega = 4$ and $\Omega = \alpha_1 + \alpha_2 + \alpha_3 + \alpha_4$.
The two theorems and the conjecture give
\begin{align*}
		\kappa_0(n) &= \frac{1}{2} \sum_{i=0}^\infty \frac{1}{2^i} 
		\binom{\alpha_1 + i}{\alpha_1} \binom{\alpha_2 + i}{\alpha_2} \binom{\alpha_3 + i}{\alpha_3} \binom{\alpha_4 + i}{\alpha_4} 	\\
		\kappa_0(n) &= 
		\sum_{i=0}^{\Omega} 
		\sum_{j=0}^i (-1)^{i - j} \binom{i}{j} 	
		\binom{\alpha_1 + j}{\alpha_1} \binom{\alpha_2 + j}{\alpha_2} \binom{\alpha_3 + j}{\alpha_3} \binom{\alpha_4 + j}{\alpha_4}	\\
		\kappa_0(n) &= 2^{\alpha_4} \! 
		\sum_{i_1=0}^{\alpha_1} 
		\sum_{i_2=0}^{\alpha_2} 
		\sum_{i_3=0}^{\alpha_3}	\!
		\binom{\! \alpha_1 \!}{i_1} 
		\binom{\! \alpha_2 \!}{i_2} 
		\binom{\! \alpha_3 \!}{i_3}
		\binom{\! \alpha_2 + i_1 \!}{\alpha_2} 
		\binom{\! \alpha_3 + i_1 + i_2 \!}{\alpha_3}	
		\binom{\! \alpha_4 + i_1 + i_2 + i_3 \!}{\alpha_4}.
\end{align*}
\begin{table}[b!]
\setlength{\tabcolsep}{1pt}
\begin{small}
\begin{tabular*}{\textwidth}{@{\extracolsep{\fill}}llcrrrrrrrrrrrrrrrrrrrr}
					&& \emph{Dirichlet}												 \\ 
					&& \emph{series}			&& \multicolumn{13}{c}{\emph{Values from} $n=1$ \emph{to} $n=12$} 			\\  
$K$					&& $1/(2 - \zeta(s))$ 			&& 1 & 1 & 1 & 2 & 1 & 3 & 1 & 4 & 2 & 3 & 1 & 8 		&& A074206			\\
$\kappa_0$			&& $\zeta(s)/(2 - \zeta(s))$ 	&& 1 & 2 & 2 & 4 & 2 & 6 & 2 & 8 & 4 & 6 & 2 & 16 		&& A067824			\\
$\tau_1\equiv \1$ 		&& $\zeta(s)$ 				&& 1 & 1 & 1 & 1 & 1 & 1 & 1 & 1 & 1 & 1 & 1 & 1 		&&					\\
$\tau_2$ 				&& $\zeta^2(s)$ 			&& 1 & 2 & 2 & 3 & 2 & 4 & 2 & 4 & 3 & 4 & 2 & 6 		&& A000005 			\\
$\tau_3$ 				&& $\zeta^3(s)$ 			&& 1 & 3 & 3 & 6 & 3 & 9 & 3 & 10 & 6 & 9 & 3 & 18 		&& A007425  			\\
$\tau_4$ 				&& $\zeta^4(s)$ 			&& 1 & 4 & 4 & 10 & 4 & 16 & 4 & 20 & 10 & 16 & 4 & 40 	&& A007426 			\\
$\upsilon_1 \equiv \1$	&& $\zeta(s)$  				&& 1 & 1 & 1 & 1 & 1 & 1 & 1 & 1 & 1 & 1 & 1 & 1 		&&					\\
$\upsilon_2$			&& $\zeta(s) (\zeta(s) - 1)$	&& 0 & 1 & 1 & 2 & 1 & 3 & 1 & 3 & 2 & 3 & 1 & 5 		&& A032741 			\\ 
$\upsilon_3$			&& $\zeta(s) (\zeta(s) - 1)^2$	&& 0 & 0 & 0 & 1 & 0 & 2 & 0 & 3 & 1 & 2 & 0 & 7 		&& A343879 			\\
$\upsilon_4$			&& $\zeta(s) (\zeta(s) - 1)^3$	&& 0 & 0 & 0 & 0 & 0 & 0 & 0 & 1 & 0 & 0 & 0 & 3 		&&					
\end{tabular*}
\end{small}
\caption{\small 
{\bf Arithmetic functions.}
For each of the arithmetic functions used in this paper, we give its Dirichlet series, the first 12 terms of its sequence, and its OEIS reference \cite{Sloane}.
Notice that $\kappa_0(n)$ is the sum of the of the $\upsilon_i$ and $\kappa_0(n) = 2 K(n)$ for $n \geq 2$.
}
\label{avalstable}
\end{table}
\section{Proof of Theorem 1}
\noindent
It is convenient to rewrite (\ref{kappadef}) as
	\begin{equation*}
		2 \kappa_0(n) = 1 + \sum_{d | n} \kappa_0(d),
	\end{equation*}
where now the sum is over all the divisors of $n$ rather than just the proper divisors.
We can express this in the language of Dirichlet convolutions:
	\begin{equation*}
		\kappa_0 = (\1 + \1 \star \kappa_0)/2,
	\end{equation*}
where $\1$ is the all 1s sequence, $1,1,1,\ldots$.
Iterating this recursive identity leads to the infinite series
	\begin{align}
		\kappa_0			&= \frac{\1}{2} + \frac{\1 \star \1}{2^2} + \frac{\1 \star \1 \star \1}{2^3} + \ldots. 
		\label{DM}	
	\end{align}
We can rewrite this as 
	\begin{align}
		\kappa_0(n)		&= \frac{\tau_1(n)}{2} + \frac{\tau_2(n)}{2^2} + \frac{\tau_3(n)}{2^3} + \ldots,
		\label{DO}	
	\end{align}
where $\tau_1 = 1,1,1,\ldots$ and 
	\begin{align*}
		\tau_i(n) 	&= \sum_{d|n} \tau_{i-1}(d). 
	\end{align*} 
Values of $\tau_1$ to $\tau_4$ are shown in Table 1. 
These quantities have a natural interpretation:
$\tau_2 \equiv d$ is the number of divisors of $n$;  
$\tau_3$ is the number of divisors of the divisors of $n$;
and so on.
It is well known that
	\begin{align*}
		\tau_2 \equiv d = \prod_{k=1}^\omega (\alpha_k + 1),
	\end{align*} 
and in general
	\begin{align}
		\tau_i 	&= \prod_{k=1}^\omega \binom{\alpha_k + i-1}{i-1}.
		\label{DQ}
	\end{align} 
Substituting this into (\ref{DO}), we obtain the desired result, namely,
	\begin{align*}
		\kappa_0(n)	&= \sum_{i=1}^\infty \frac{1}{2^i} \prod_{k=1}^\omega \binom{\alpha_k + i -1}{i-1} \\
					&= \frac{1}{2} \sum_{i=0}^\infty \frac{1}{2^i} \prod_{k=1}^\omega \binom{\alpha_k + i}{i}.
	\end{align*}
This can be expressed as a generalized hypergeometric function:
	\begin{equation*}
		\kappa_0(n) = \textstyle{\frac{1}{2}} \, {}_\omega F_{\omega-1}\!
		\left[
			\begin{matrix}
				\alpha_1+1, \ldots, \alpha_\omega+1\\1, \ldots, 1
			\end{matrix}
			; \frac{1}{2}
		\right]. 
		\qed
	\end{equation*}
\section{Proof of Theorem 2}
\noindent
We know that $\kappa_0(n)$ is the total number of squares in the divisor tree of $n$ (Fig. 1E).
Let's consider the construction of a divisor tree, generation by generation, as shown in Fig. 1.
Let $\upsilon_1(n) = 1$ be the number of squares in the root of the tree, namely, the single largest square (Fig 1A).
Let $\upsilon_2(n)$ be the number of squares in the main arm of the tree, not including the root (Fig 1B).
Let $\upsilon_3(n)$ be the number of squares in the secondary arms, not including their roots (Fig. 1C), and so on.
The quantity $\upsilon_i(n)$ has a natural interpretation:
$\upsilon_2(n)$ is the number of proper divisors of $n$,
$\upsilon_3(n)$ is the number of proper divisors of the proper divisors of $n$, and so on.
As with $\tau_i(n)$, we have $\upsilon_1 = 1,1,1,\ldots$ and 
	\begin{align}
		\upsilon_i(n) &= \sum_{d\divprop n} \upsilon_{i-1}(d).
		\label{FR}
	\end{align} 
Values of $\upsilon_1$ to $\upsilon_4$ are shown in Table 1. 
We can then express $\kappa_0(n)$ as the sum of the $\upsilon_i$ over the root and the $\Omega(n)$ generations of arms:
	\begin{align}
		\kappa_0(n) = \sum_{i=0}^{\Omega} \upsilon_{i+1}.
		\label{GA}
	\end{align}
The $\upsilon_k$ are related to the $\tau_k$ as follows:
	\begin{align*}
		\upsilon_1 &= \tau_1 \\
		\upsilon_2 &= \tau_2 - \tau_1 \\
		\upsilon_3 &= \tau_3 - 2 \tau_2 + \tau_1\\
		\upsilon_4 &= \tau_4 - 3 \tau_3 + 3 \tau_2 - \tau_1,
	\end{align*}
and in general, by the principle of inclusion and exclusion,
\begin{align*}
\upsilon_i &= \sum_{j=0}^{i-1} (-1)^{i-1-j} \binom{i-1}{j} \tau_{j+1}.
\end{align*}
Substituting this into (\ref{GA}), 
	\begin{align*}
		\kappa_0(n) &= \sum_{i=0}^{\Omega} \sum_{j=0}^{i} (-1)^{i-j} \binom{i}{j} \tau_{j+1}.
	\end{align*}
Substituting $\tau_i$ from (\ref{DQ}) into this, we obtain the desired result:
\begin{align*}
\kappa_0(n)	&= \sum_{i=0}^{\Omega} \sum_{j=0}^{i} (-1)^{i-j} \binom{i}{j} \prod_{k=1}^\omega \binom{\alpha_k + j}{j}. \qed
\end{align*}
\section{Discussion}
\noindent
Theorem 1 tells us that, to our surprise, $\kappa_0(n) = 2 K(n)$ can be expressed as a generalized hypergeometric function.
It is in some sense the simplest such function that can be naturally tied to the prime signature of a number.
This connection opens the door to the considerable machinery that is known for the generalized hypergeometric function.
\\ \indent
Theorem 1 offers some insight into the properties of $\kappa_0(n)$ and $K(n)$.
Our proof of it suggests a simple demonstration that, for $n \geq 2$, $\kappa_0(n) = 2 K(n)$. 
As with $\kappa_0(n)$,  it is convenient to express (\ref{kappadef}) as
	\begin{equation*}
		2 K(n) = \varepsilon + \sum_{d | n} K(d).
	\end{equation*}
In the language of Dirichlet convolutions, this is
	\begin{equation}
		K = (\varepsilon + \1 \star K)/2,
		\label{GG}
	\end{equation}
where recall $\varepsilon = 1,0,0,0,\ldots$. 
Iterating this leads to the infinite series
	\begin{align*}
		K 	&= \frac{\varepsilon}{2} + \frac{\1}{2^2} + \frac{\1 \star \1}{2^3} + \frac{\1 \star \1 \star \1}{2^4} + \ldots \\
		 	&= \frac{1}{2} \left(\varepsilon + \frac{\1}{2} + \frac{\1 \star \1}{2^2} + \frac{\1 \star \1 \star \1}{2^3} + \ldots \right) \\
			&= \left(\varepsilon + \kappa_0 \right)/2,
	\end{align*}
where the last step makes use of (\ref{DM}).
\\ \indent
We can also readily calculate the Dirichlet series for $\kappa_0(n)$ and $K(n)$.
Our starting point is  (\ref{DO}).
Denoting the Dirichlet series of $\kappa_0$ and $K$ by $\widetilde \kappa_0$ and $\widetilde K$, we can write
	\begin{align*}
		\widetilde \kappa_0	&= \frac{\widetilde \tau_1}{2} + \frac{\widetilde \tau_2}{2^2} + \frac{\widetilde \tau_3}{2^3} + \ldots.
	\end{align*}
Since $\tau_1 = \1$, $\tau_2 = \1 \star \1$, and so on,
and the Dirichlet series for $\1$ is $\zeta(s)$, we have $\widetilde \tau_i = \zeta^i(s)$.
Then
	\begin{align*}
		\widetilde \kappa_0	&= \frac{\zeta(s)}{2} + \frac{\zeta(s)^2}{2^2} + \frac{\zeta(s)^3}{2^3} + \ldots \\
						&= \frac{\zeta(s)}{2 - \zeta(s)}.
	\end{align*}
As for $\widetilde K$, from (\ref{GG}), $2 \widetilde K = 1 + \widetilde \kappa_0$, so $\widetilde K = 1/(2 - \zeta(s))$.
\\ \indent
When $n = p_1 p_2 \ldots p_\omega$ is the product of $\omega$ distinct primes, 
all of the $\alpha_i$ equal one.
Then Theorem 1 reduces to 
	\begin{align*}
		\kappa_0(p_1 \ldots p_\omega) 	&= \frac{1^\omega}{2} + \frac{2^\omega}{2^2} + \frac{3^\omega}{2^3} + \ldots \\
						&= {\rm Li}_{-\omega}(1/2),
	\end{align*}
where {\rm Li} is the polylogarithm.
For $\omega=1, 2, 3, \ldots$, this has values 2, 6, 26, 150, 1082, \ldots.
Its exponential generating function is
	\begin{align*}
		{\rm EG}(\kappa_0(p_1 \ldots p_\omega),x) &= \sum_{\omega = 0}^\infty \frac{x^\omega}{\omega!} \sum_{i=1}^\infty \frac{i^\omega}{2^i} \\
				&= \sum_{i = 1}^\infty \frac{1}{2^i} \sum_{\omega=0}^\infty \frac{(i x)^\omega}{\omega!}	\\
				&= \sum_{i = 1}^\infty \frac{e^{ix}}{2^i} \\
				&= \frac{e^x}{2 - e^x}.
	\end{align*}
\indent
MacMahon \cite{MacMahon93} derived a somewhat more complex version of Theorem 2,
using a more laborious approach. 
It is 
	\begin{align*}
		K(n) & = \sum_{i=1}^{\Omega} \sum_{j=0}^{i-1} (-1)^{j} \binom{i}{j} 	
		\prod_{k=1}^\omega \binom{\alpha_k + i - j - 1}{\alpha_k}.
	\end{align*}
\\ \indent
In proving Theorem 2, we made use of the sequences $\upsilon_i$, shown in Table 1.
We can also calculate their Dirichlet series.
Adding $\upsilon_{i-1}$ to (\ref{FR}), and turning to the language of Dirichlet convolutions,
	\begin{align*}
		\upsilon_i + \upsilon_{i-1} &= \1 \star \upsilon_{i-1}.
	\end{align*}
Denoting the Dirichlet series of $\upsilon_i$ by $\widetilde \upsilon_i$, this implies
$\widetilde \upsilon_i + \widetilde \upsilon_{i-1} = \zeta(s) \widetilde\upsilon_{i-1}$, that is,
	\begin{align*}
		\widetilde \upsilon_i 	&= (\zeta(s) - 1) \widetilde \upsilon_{i-1}.
	\end{align*}
Since $\widetilde \upsilon_1 = \zeta(s)$, it follows that $$\widetilde \upsilon_i = \zeta(s) (\zeta(s) - 1)^{i-1}.$$
\\ \indent
Conjecture 1, which is correct for $n \leq 10^5$, can be written more symmetrically:
	\begin{align*}
		\kappa_0(n) &=	\sum_{i_1=0}^{\alpha_1} 	\ldots 	\sum_{i_{\omega}=0}^{\alpha_{\omega}} 
		\prod_{k=1}^{\omega}	
		\binom{\alpha_k}{i_k} \binom{\alpha_k + i_1 + \ldots + i_{k-1}}{\alpha_k}.	
	\end{align*}
But since $i_\omega$ only appears in $\binom{\alpha_\omega}{i_\omega}$, it sums to $2^{\alpha_\omega}$, giving the original form of Conjecture 1.
Since the $\alpha_k$ can be permuted at will, a corollary of this is that $\kappa_0(n)$ is divisible by $2^{\alpha^*}$,
where $\alpha^*$ is the largest of the $\alpha_k$s, which we proved in \cite{Fink21}.
This makes Conjecture 1 the most efficient of our three expressions for calculating values of $K$ and $\kappa_0$ for very large values of $n$. 
In particular, it is useful for calculating the indices of the sequence records of $K$ and $\kappa_0$---the K-champion numbers \cite{Deleglise08} (A307866 \cite{Sloane}) and the recursively highly composite numbers \cite{Fink21} (A333952 \cite{Sloane}).
\\ \indent
There are a number of open questions about $K(n)$ and $\kappa_0(n)$.
Here are four.
\\ 1. Can Conjecture 1 be proved?
\\ 2. Can Theorem 1 be be generalized from $\kappa_0(n)$ to $\kappa_x(n)$?
\\ 3. When $\alpha_i = 1$ for all $i$, 
Theorem 1 reduces to the polylogarithm and has a simple exponential generating function.
What about when $\alpha_i = j$ for all $i$?
\\ 4. What is the significance of the sequence $\kappa_0(n)/2^{\alpha^*(n)}$, namely,
1, 1, 1, 1, 1, 3, 1, 1, 1, 3, 1, 4, 1, 3, 3, 1, 1, 4, 1, 4, 3, 3, 1, 5, 1, 3, 1, 4, 1, 13?
(The last and 30th term distinguishes this sequence from others.)



\begin{thebibliography}{1}
\begin{footnotesize}
\bibitem{Kalmar31}		L. Kalmar, 							A factorisatio numerorum probelmajarol, 										\emph{Mat Fiz Lapok}		{\bf 38}, 	1 		(1931).
\bibitem{Hille36} 		E. Hille,								A problem in factorisatio numerorum, 										\emph{Acta Arith} 			{\bf 2}, 	134		(1936).
\bibitem{Canfield83} 		E. Canfield, P. Erd\"{o}s, C, Pomerance,		On a problem of Oppenheim concerning ``factorisatio numerorum",					\emph{J Number Theory} 		{\bf 17},  	1	 	(1983).
\bibitem{Deleglise08}	M. Del\'{e}glise, M.\ Hernane, J.-L. Nicolas,	Grandes valeurs et nombres champions de la fonction arithm\'{e}tique de Kalm\'{a}r, 		\emph{J Number Theory} 		{\bf 128}, 	1676 	(2008).
\bibitem{Newberg93} 	L. Newberg, D. Naor, 					A lower bound on the number of solutions to the probed partial digest problem,			\emph{Adv Appl Math} 		{\bf 14},	172		(1993).
\bibitem{Chor00} 		B. Chor, P. Lemke, Z. Mador,				On the number of ordered factorizations of natural numbers, 						\emph{Disc Math} 			{\bf 214},	123		(2000).
\bibitem{Klazar07} 		M. Klazar, F. Luca,						On the maximal order of numbers in the factorisatio numerorum problem,				\emph{J Number Theory} 		{\bf 124},  470	 	(2007).
\bibitem{Fink21} 		T. Fink,								Recursively divisible numbers,												arxiv.org/abs/1912.07979.
\bibitem{Fink22} 		T. Fink,								Recursively abundant and recursively perfect numbers,							arxiv.org/abs/2008.10398.
\bibitem{Sloane}		N.\ J.\ A.\ Sloane, editor, 					The On-Line Encyclopedia of Integer Sequences, published electronically at https://oeis.org, 2018.
\bibitem{MacMahon93} 	P. A. MacMahon,						Memoir on the theory of the compositions of numbers,							\emph{Philos T R Soc Lond A} 	{\bf 184}, 	835 		(1893).
\end{footnotesize}
\end{thebibliography}
\end{document}